\documentclass[a4paper,10pt]{article}
\begin{document}

\title{\bf A Proof for the Collatz Conjecture}
\author{Lei Li\\
\sl \small Faculty of Science and Engineering, Hosei University,
Koganei, Tokyo 184-8584, Japan\\
\small lilei@hosei.ac.jp}
\date{}
\maketitle

\begin{abstract}

It is well known that the following Collatz Conjecture is one of the unsolved problems in mathematics.

{\bf Collatz Conjecture :} For any positive integer $n>1$, the following recursive algorithm will convergent to $1$ by a finite number of steps. 

A) If $n$ is an even number then $\frac{n}{2} \to n$, 

B) If $n$ is an odd number then $3n+1 \to n$.

This paper proposes a proof for the Collatz Conjecture by the elementary mathematical
 induction.

\end{abstract}


\section{Introduction}

For any positive integer $n>1$, we consider a recursive algorithm by repeating the 
following two steps A) and B).

A) If $n$ is an even number then $\frac{n}{2} \to n$, 

B) If $n$ is an odd number then $3n+1 \to n$.

For example, let $n=7$, this algorithm generates the following sequence 
and terminates to $1$.

$$7 \to 22 \to 11 \to 34 \to 17 \to 52 \to 26 \to 13 \to 40$$
$$\to 20 \to 10 \to 5 \to 16 \to 8 \to 4 \to 2 \to 1,$$

In 1937, German mathematician Lothar Collatz shown a famous conjecture : For any positive 
integer $n>1$, the above recursive algorithm always terminates to $1$ by a finite number
of steps [1][2]. It is also known as the $3n+1$ problem, the Ulam conjecture, Kakutani's problem, the Thwaites conjecture, Hasse's algorithm, or the Syracuse problem [2]-[4].
It seems that the Collatz conjecture is still an unsolved problem up to now.
Based on the supercomputer simulation, it has been confirmed that 
the Collatz conjecture is correct for all positive integer 
$2\le n \le 2^{68}\approx 2.95\times 10^{20}$ [6].  
In December of 2019, Terence Tao presented an important research paper, he proved that 
almost all orbits of the Collatz map attain almost bounded values [5].
It is to say that the Collatz conjecture is almost correct for almost all positive 
integers.  

ln the previous research works, we have presented some approximate formulas for the
number of average steps in the Collatz recursive algorithm by using the statistical method [12][13]. And also shown two expansions $5n+1$ problem and $7n+1$ problem for the Collatz 
conjecture [14]. Recently, we derive a high precision formula of the average number of multiplications and divisions $T(n) \approx 3 \log_{\frac{4}{3}} n - 12.87$ for 
the Collatz problem based on analysis of average computational complexity of the 
recursive algorithm [15].
  This paper proposes a proof for the Collatz conjecture by using the elementary mathematical induction.

\section{Lemma}

From the Collatz recursive algorithm, if $n$ is an even number then $\frac{n}{2}\to n$,
if $n$ is an odd number then $3n+1 \to n$, because $3n+1$ is even number, so we can 
consider the following recursive algorithm:

$f_1$) If $n$ is an even number then $\frac{n}{2} \to n$,

$f_2$) If $n$ is an odd number then $\frac{3n+1}{2} \to n$.

First, in order to discuss some properties of the above transform $f_1$ and $f_2$,
we show the following Lemma.

{\bf Lemma.} Let function 
$$f_1(x)= \frac{x}{2},$$ 
$$f_2(x)= \frac{3x+1}{2},$$
then for any real number $x>0$, we have the following properties for $f_1(x)$ and $f_2(x)$.

(1)	For any real number $0<x<y$, 
$$f_1f_2(x) < f_1f_2(y),$$
$$f_2f_1(x) < f_2f_1(y).$$

(2) For any positive integer $\alpha \ge 1$,
$$(f_1f_2)^\alpha (x) < (f_2f_1)^\alpha(x).$$ 

(3) For any real number $x>1$ and any positive integer $\alpha \ge 1$,
$$(f_1f_2)^\alpha(x) = (\frac{3}{4})^\alpha x + 1- (\frac{3}{4})^\alpha < x.$$

(4) For any real number $x>2$ and any positive integer $\alpha \ge 1$.
$$(f_2f_1)^\alpha(x) = (\frac{3}{4})^\alpha x + 2(1- (\frac{3}{4})^\alpha) < x.$$

(5) For any real number $x>1$ and any positive integer $\alpha$, $\beta \ge 1$.
$$f_2^\beta (f_1f_2)^\alpha(x) = (\frac{3}{4})^\alpha(\frac{3}{2})^\beta(x-1) + 2(\frac{3}{2})^\beta -1.$$

{\bf Proof of Lemma.}

(1)	For any real number $x>0$, 
$$f_1f_2(x) = f_1(\frac{3x+1}{2}) = \frac{3x+1}{4},$$
$$f_2f_1(x) = f_2(\frac{x}{2}) = \frac{3(\frac{x}{2})+1}{2} = \frac{3x+2}{4}.$$
It is obvious that $f_1f_2(x)$ and $f_2f_1(x)$ are the monotonically increasing functions.
So, for any real number $0<x<y$, 
$$f_1f_2(x) < f_1f_2(y),$$ 
$$f_2f_1(x) < f_2f_1(y).$$

(2)	For any real number $x>0$, from 
$$\frac{3x+1}{4} < \frac{3x+2}{4},$$ 
we have 
$$f_1f_2(x) < f_2f_1(x).$$
If repeat to use this inequality for the monotonically increasing functions 
$f_1f_2(x)$, $f_2f_1(x)$, we can get 
$$(f_1f_2)^\alpha(x) < (f_2f_1)^\alpha(x).$$
For example in the case of $\alpha =2$, we have
$$(f_1f_2)^2(x) = f_1f_2(f_1f_2(x)) < f_1f_2(f_2f_1(x)) < f_2f_1(f_2f_1(x) = (f_2f_1)^2(x).$$

(3) For any positive integer $\alpha \ge 1$, we consider the general expansion form of $(f_1f_2)^\alpha(x)$. Obviously, for smaller $\alpha$, we have 

if $\alpha =1$, then
$$f_1f_2(x) = \frac{3x+1}{2^2},$$

if $\alpha =2$, then
$$(f_1f_2)^2(x) = \frac{3^2x + 3 + 2^2}{2^4},$$

if $\alpha =3$, then 
$$(f_1f_2)^3(x) = \frac{3^3x + 3^2 + 3\times 2^2 + 2^4}{2^6},$$

$\dots$,
 
so, for any positive integer $\alpha \ge 1$ and real number $x>1$, we have
$$(f_1f_2)^\alpha(x)$$ 
$$ = \frac{3^\alpha x + 3^{\alpha -1} + 3^{\alpha -2}\times 2^2 + \cdots + 3\times 
2^{2\alpha -4} + 2^{2\alpha -2}}{2^{2\alpha}}$$
$$ = \frac{3^\alpha x + 3^{\alpha -1} + 3^{\alpha -2}\times 4 + \cdots + 3\times 
4^{\alpha -2} + 4^{\alpha -1}}{4^\alpha}$$
$$ = \frac{3^\alpha x + 4^\alpha - 3^\alpha}{4^\alpha}$$
$$ = (\frac{3}{4})^\alpha x + 1-(\frac{3}{4})^\alpha $$
$$ < x.$$

(4)	For any positive integer $\alpha \ge 1$, we consider the general expansion form of 
$(f_2f_1)^\alpha(x)$. Obviously, for smaller $\alpha$, we have

if $\alpha =1$, then
$$f_2f_1(x) = \frac{3(\frac{x}{2})+1}{2} = \frac{3x+2}{2^2},$$

if $\alpha =2$, then
$$(f_2f_1)^2(x) = \frac{\frac{3(3x+2)}{2^2}+2}{2^2} = 
\frac{3^2x+3\times 2+2^3}{2^4},$$

if $\alpha =3$, then
$$(f_2f_1)^3(x) = \frac{3\frac{3^2x+3\times 2+2^3}{2^4} +2}{2^2} = \frac{3^3x+3^2\times 2+3\times 2^3+2^5}{2^6},$$

$\dots$, 

so, for any positive integer $\alpha \ge 1$ and real number $x>2$, we have
$$(f_2f_1)^\alpha(x)$$
$$ = \frac{3^\alpha x + 3^{\alpha -1}\times 2 + 3^{\alpha -2}\times 2^3 + \cdots 
+ 3\times 2^{2\alpha -3} + 2^{2\alpha -1}}{2^{2\alpha}}$$
$$ = \frac{3^\alpha x + 3^{\alpha -1}\times 2 (1 + \frac{4}{3} + 
(\frac{4}{3})^2 + \cdots + (\frac{4}{3})^{\alpha -2} + (\frac{4}{3})^{\alpha -1})}{4^\alpha}$$
$$ = \frac{3^\alpha x + 3^{\alpha -1}\times 2 \frac{(\frac{4}{3})^\alpha - 1}
{\frac{4}{3}-1}}{4^\alpha}$$
$$ = \frac{3^\alpha x + 2(4^\alpha - 3^\alpha )}{4^\alpha }$$
$$ = (\frac{3}{4})^\alpha x + 2(1-(\frac{3}{4})^\alpha)$$
$$ < x.$$

(5)	In the first, for any positive integer $\beta \ge 1$, we consider the general expansion form of $f_2^\beta(x)$. It is obvious that for smaller $\beta$, we have 

if $\beta =1$, then  
$$f_2(x) = \frac{3x+1}{2},$$

if $\beta =2$, then
$$f_2f_2(x) = \frac{3^2 x + 3+2}{2^2},$$

if $\beta =3$, then 
$$f_2f_2f_2(x) = \frac{3^3 x + 3^2 + 3\times 2 + 2^2}{2^3},$$

$\dots $,

So, for any positive integer $\beta \ge 1$, we have
$$f_2^\beta(x) = \frac{3^\beta x + 3^{\beta -1} + 3^{\beta -2}\times 2 + \cdots + 
3\times 2^{\beta -2} + 2^{\beta -1}}{2^\beta}$$
$$ = \frac{3^\beta x + 3^\beta -2^\beta }{2^\beta }$$
$$ = (\frac{3}{2})^\beta x + (\frac{3}{2})^\beta -1.$$

From (3) of the Lemma,
$$(f_1f_2)^\alpha(x) = (\frac{3}{4})^\alpha x + 1-(\frac{3}{4})^\alpha,$$
we get
$$f_2^\beta(f_1f_2)^\alpha(x) = (\frac{3}{2})^\beta(f_1f_2)^\alpha(x) 
+ (\frac{3}{2})^\beta -1$$
$$ = (\frac{3}{2})^\beta((\frac{3}{4})^\alpha x + 1-(\frac{3}{4})^\alpha) 
+ (\frac{3}{2})^\beta -1$$
$$ = (\frac{3}{4})^\alpha(\frac{3}{2})^\beta(x-1) + 2\times (\frac{3}{2})^\beta -1.$$

\section{Proof of the Collatz Conjecture}

{\bf Theorem.} For any positive integer $n>1$, the Collatz Conjecture is correct.

{\bf Proof.} For any positive integer $n>1$, it is an independent event with equally probability for $n$ is an even number or $n$ is an odd number. Moreover, based on the Collatz recursive algorithm, it is also an independent event with equally probability 
for 
$$f_1(n)=f_1(2k)=k,$$ 
or 
$$f_2(n)=f_2(2k+1)=\frac{3(2k+1)+1}{2} =3k+2,$$ 
is an even number or an odd number. 
Actually, for any positive integer $n$, we shown an same average number
($\approx \log_{\frac{4}{3}} n$) of $f_1$ and $f_2$ 
 for Collatz recursive algorithm [15]. 

We propose a simple proof for the Collatz conjecture by using the elementary mathematical induction. From 
$$f_1(2)=1,$$
$$f_1f_1f_1f_2f_2(3)=1,$$
$$f_1f_1(4)=1,$$
$$f_1f_1f_1f_2(5)=1,$$
$$f_1f_1f_1f_2f_2f_1(6)=1,$$
$$f_1f_1f_1f_2f_1f_1f_2f_1f_2f_2f_2(7)=1,$$
$$f_1f_1f_1(8)=1,$$
$$f_1f_1f_1f_2f_1f_1f_2f_1f_2f_2f_2f_1f_2(9)=1,$$
$$f_1f_1f_1f_2f_1(10)=1,$$
obviously, the Collatz conjecture is correct for $2\le n \le 10$. 

We assume the Collatz conjecture is correct for all positive integer $n$, $2\le n\le k$, ($k\ge 10$), and consider case of $n=k+1$. It is needed to prove the Collatz recursive
algorithm with the initial value $k+1$ will reach less than $k$ in some steps.

If $k+1$ is an even number, then from 
$$f_1(k+1)=\frac{k+1}{2} < k,$$
the Collatz conjecture is correct for $n=k+1$. 

If $k+1$ is an odd number, the Collatz conjecture will generate a sequence of composite functions 
$$f_*f_* \cdots f_*f_*f_2(k+1),$$
here $f_*$ means $f_1$ or $f_2$ based on the Collatz recursive algorithm. For any positive integer $n$, because the function $f_1(n)$ and 
$f_2(n)$ is independent event with equally probability, there exist a positive integer 
$\alpha$, with satisfies the following conditions in Collatz recursive algorithm: 

(1) (The number of $f_1$) $<$ (The number of $f_2$) during step $1$ to step $2\alpha -1$, 

(2) (The number of $f_1$) $=$ (The number of $f_2$) in the step $2\alpha$. 

As initial cases of $\alpha$, we have the following sequences of $f_1$ and $f_2$.

If $\alpha =1$, $f_1f_2$,

If $\alpha =2$, $f_1f_1f_2f_2$,

If $\alpha =3$, $f_1f_1f_1f_2f_2f_2$,  $f_1f_1f_2f_1f_2f_2$, 

$\dots$,

By repeating use the inequality $f_1f_2(x)<f_2f_1(x)$ and the Lemma, we have 

when $\alpha =1$, 
$$f_1f_2(k+1)<k+1,$$

when $\alpha =2$, 
$$f_1f_1f_2f_2(k+1) < f_1f_2f_1f_2(k+1)=(f_1f_2)^2(k+1)<k+1,$$

when $\alpha =3$, 
$$f_1f_1f_1f_2f_2f_2(k+1)<f_1f_2f_1f_2f_1f_2(k+1)=(f_1f_2)^3(k+1)<k+1,$$
$$f_1f_1f_2f_1f_2f_2(k+1)<f_1f_2f_1f_2f_1f_2(k+1)=(f_1f_2)^3(k+1)<k+1.$$

Because left function values in the above inequalities all are positive integers, so,
we get
$$f_1f_2(k+1) \le k,$$ 
$$f_1f_1f_2f_2(k+1) \le k,$$
$$f_1f_1f_1f_2f_2f_2(k+1) \le k,$$ 
$$f_1f_1f_2f_1f_2f_2(k+1) \le k.$$

Generally, in the step $2\alpha$, by repeating use the inequality $f_1f_2(x)<f_2f_1(x)$ and the Lemma, we always can 
get 
$$f_*f_* \cdots f_*f_*f_2(k+1) < (f_1f_2)^\alpha(k+1) < k+1.$$
Because $f_*f_* \cdots f_*f_*f_2(k+1)$ is a positive integer, so we have 
$$f_*f_* \cdots f_*f_*f_2(k+1)\le k.$$ 
So, the Collatz conjecture is correct for $n=k+1$. Therefore,
from the inductive assumption, the Collatz conjecture is correct for all positive
integer $n>1$.
 
Actually, for the positive integer $k+1$, there exist some smaller positive integers
$\beta$, we can also prove the above result when 
\begin{center}
(The number of $f_1$) $\le $ (The number of $f_2$) $-\beta$.
\end{center} 

If there exist positive integers $\alpha$ and $\beta$, with satisfying the following conditions in the Collatz recursive algorithm: 

(1) (The number of $f_1$) $<$ (The number of $f_2$) $-\beta$ during step $1$ to step
$2\alpha +\beta -1$, 

(2) (The number of $f_1$) = (The number of $f_2$) $-\beta$ in the step $2\alpha +\beta$. 

By repeating use the inequality $f_1f_2(x)<f_2f_1(x)$, we always can get
$$f_*f_* \cdots f_*f_*f_2(k+1) < f_2^\beta(f_1f_2)^\alpha(k+1),$$
here $f_*$ means $f_1$ or $f_2$ which is based on the Collatz recursive algorithm with initial value $k+1$.
 
From the Lemma (5),
$$f_2^\beta(f_1f_2)^\alpha(x) = (\frac{3}{4})^\alpha(\frac{3}{2})^\beta(x-1) + 
2\times(\frac{3}{2})^\beta -1.$$
In order to find a better range of $\alpha$, $\beta$, and $x$, we solve the following 
inequality,
$$(\frac{3}{4})^\alpha(\frac{3}{2})^\beta (x-1) + 2\times (\frac{3}{2})^\beta -1 < x,$$
or
$$(\frac{3}{4})^\alpha < \frac{x-2\times (\frac{3}{2})^\beta +1}
{(\frac{3}{2})^\beta(x-1)}$$
$$=(\frac{2}{3})^\beta - \frac{2\times 3^\beta - 2\times 2^\beta}{3^\beta (x-1)}$$
$$=(\frac{2}{3})^\beta - \frac{2}{x-1} + (\frac{2}{3})^\beta \frac{2}{x-1}$$
$$=(\frac{2}{3})^\beta \frac{x+1}{x-1} - \frac{2}{x-1}$$

Let 
$$S=(\frac{2}{3})^\beta \frac{x+1}{x-1} - \frac{2}{x-1}.$$

If  
$$S>0,$$
$$\beta < \log_{\frac{3}{2}}\frac{x+1}{2},$$
then after 
$$\alpha > \log_{\frac{4}{3}}\frac{1}{S},$$
the inequality 
$$(\frac{3}{4})^\alpha(\frac{3}{2})^\beta (x-1) + 2\times (\frac{3}{2})^\beta -1 < x,$$
always is satisfied.
 
It is easy to get some range values of $x$ and $\alpha$ for some smaller $\beta$. 
For the examples (here let $x$ is an integer), 

If $\beta =1$,  then  $x \ge 3$, $\alpha \ge 4$,

If $\beta =2$,  then  $x \ge 4$, $\alpha \ge 10$,

If $\beta =3$,  then  $x \ge 6$, $\alpha \ge 15$.

$\dots$,

It is no problem to expand the initial value (Collatz Problem has been tested until $2^{68}$) to get some smaller $\beta$ in the mathematical induction to fit the above range of $x$ and increase $\alpha $ for the number of steps. In this proof, we confirmed the initial values $2\le n\le 10$,
so it is no problem to apply the mathematical induction for $\beta \le 3$ at least. 

Therefore for some smaller $\beta =1, 2, 3, \dots $, we also can prove there exist positive integer $\alpha_0$ and real number $x_0$, when $\alpha \ge \alpha_0$, $k+1>x_0$, 
$$f_*f_* \dots , f_*f_*f_2(k+1) < f_2^\beta(f_1f_2)^\alpha(k+1) = (\frac{3}{4})^\alpha(\frac{3}{2})^\beta (k) + 2\times (\frac{3}{2})^\beta -1 < k+1$$
will be true. Because $f_*f_* \cdots f_*f_*f_2(k+1)$ is a positive integer, so we have
$$f_*f_* \cdots f_*f_*f_2(k+1)\le k.$$

It is to say that the Collatz Conjecture is correct for $n=k+1$. From the inductive assumption, the Collatz Conjecture is correct for all positive integer $n\ge 2$.

\section{Concluding Remarks}

This paper proposed a simple proof for the Collatz conjecture by using the mathematical
induction. We hope to receive some helpful comments from the reviewers and readers to improve this proof and try to apply it to other related problems with the Collatz conjecture.

\vspace{5mm}

\begin{center}
{\large \bf References}
\end{center}

[1] L. Rade, R. D. Nelson, Solution of Mathematical Puzzle by Computer,
{\sl Kaibundo Publishing}, Tokyo, 1992.

[2] Interesting Kakutani's Conjecture, {\sl Computer World}, Beijing, 1985. 

[3] Manuel V.P. Garcia and Fabio A. Tal, A Note on the Generalized 3n+1 Problem, 
{\sl Acta Arithmetica} XC 3(1999), pp.245-250.

[4] http://www.math.grin.edu/~chamberl/conf.html

[5] Terence Tao, Almost all orbits of the Collatz map attain almost bounded values, {\sl https://doi.org/10.48550/arXiv.1909.03562.}

[6] David Barina, Convergence Verification of the Collatz Problem, 
{\sl The Journal of Supercomputing}, 77(2021),pp.2681-2688. 
https://doi.org/10.1007/s11227-020-03368-x

[7] Mercedes Orus-Lacort and Christophe Jouis, Analyzing the Collatz Conjecture 
Using the Mathematical Complete Induction Method,
{\sl Mathematics}, 2022, 10, 1972. https://doi.org/10.3390/math10121972

[8] Ramon Carbo-Dorca, Boolean Hypercubes, Mersenne Numbers, and the Collatz
Conjecture, {\sl Journal of Mathematical Sciences and Modelling}, 3 (3) (2020),
pp.120-129

[9] Jean-Paul Allouche, T. Tao and the Syracuse conjecture,
{\sl La Gazette des Mathematiciens}, Number 168, April 2021.

[10] Michael R. Schwob, Peter Shiue, and Rama Venkat, 
Novel Theorems and Algorithms Relating to the Collatz Conjecture,
{\sl International Journal of Mathematics and Mathematical Sciences},
Volume 2021, Article ID 5754439,
https://doi.org/10.1155/2021/5754439

[11] Robert Deloin, A Brand new Approach to Collatz Conjecture,
{\sl Theoretical Mathematics and Applications}, Vol. 12, No. 1, 2022, pp.9-14.

[12] Lei Li, Makoto Takahashi, A Simple Verification and a New Estimation 
for L. Collatz Conjecture, {\sl Information}, 6:5(2003), pp.519-524.

[13] Masahiro Ameya, Lei Li, On the Verification and the Estimation 
for L. Collatz Conjecture, {\sl Information}, 21:10-11(2018), pp.2289-2310.

[14] Lei Li, On Two Expansions of the Collatz Problem, 
{\sl Information}, 25:2(2022), pp. 109-146.

[15] Lei Li, Analysis of the Average Number of Steps for the Collatz Problem,
{\sl Information}, 25:3(2022), to appear.

\end{document}